\documentclass[11pt]{amsart}

\usepackage[T1]{fontenc}
\usepackage{amssymb}
\usepackage{tikz}
\usetikzlibrary{arrows.meta}
\usepackage[colorlinks=true,linkcolor=blue,citecolor=blue]{hyperref}
\hypersetup{pdftitle={A right-quantum matrix-tree theorem},pdfauthor={Guo-Niu Han}}

\def\bw#1#2{\genfrac{(}{)}{0pt}{1}{#1}{#2}}
\makeatletter
\newcommand{\bwc}[1]{{\def\bw@t{}\def\bw@b{}\bw@loop#1,\relax,\bw{\bw@t}{\bw@b}}}
\def\bw@loop#1,{\ifx\relax#1\else\bw@one#1\@nil\expandafter\bw@loop\fi}
\def\bw@one#1/#2\@nil{\setbox0\hbox{$\scriptstyle#1$}\setbox2\hbox{$\scriptstyle#2$}%
  \dimen0=\wd0 \ifdim\wd2>\dimen0 \dimen0=\wd2\fi
  \ifx\bw@t\empty\else\g@addto@macro\bw@t{\;\;}\g@addto@macro\bw@b{\;\;}\fi
  \edef\bw@w{\the\dimen0}%
  \expandafter\g@addto@macro\expandafter\bw@t\expandafter{\expandafter
    \makebox\expandafter[\bw@w]{$\scriptstyle#1$}}%
  \expandafter\g@addto@macro\expandafter\bw@b\expandafter{\expandafter
    \makebox\expandafter[\bw@w]{$\scriptstyle#2$}}}
\makeatother
\def\setA{\mathbb{A}}
\def\calA{\mathcal{A}}
\def\calI{\mathcal{I}}
\def\calR{\mathcal{R}}
\def\Src{\mathop{\mathrm{Src}}\nolimits}
\def\Reg{\mathop{\mathrm{Reg}}\nolimits}
\def\Ord{\mathop{\mathrm{Ord}}\nolimits}
\def\ch{\mathop{\mathrm{ch}}\nolimits}
\def\inv{\mathop{\mathrm{inv}}\nolimits}
\def\cyc{\mathop{\mathrm{cyc}}\nolimits}
\def\rk{\mathop{\mathrm{rk}}\nolimits}
\def\Re{\mathop{\mathrm{Re}}\nolimits}
\def\B{\mathbf{B}}
\def\I{\mathbf{I}}
\def\Ferm{\mathop{\mathrm{Ferm}}\nolimits}
\def\Bos{\mathop{\mathrm{Bos}}\nolimits}

\newtheorem{theorem}{Theorem}[section]
\newtheorem{lemma}[theorem]{Lemma}
\newtheorem{proposition}[theorem]{Proposition}
\newtheorem{corollary}[theorem]{Corollary}
\theoremstyle{definition}
\newtheorem{definition}[theorem]{Definition}
\newtheorem{example}[theorem]{Example}
\newtheorem{remark}[theorem]{Remark}

\numberwithin{equation}{section}

\begin{document}

\title[A right-quantum matrix-tree theorem]
{A right-quantum matrix-tree theorem}

\author{Guo-Niu Han}
\address{I.R.M.A., UMR 7501, Universit\'e de Strasbourg et CNRS,
7 rue Ren\'e-Descartes, F-67084 Strasbourg, France}
\email{guoniu.han@unistra.fr}

\date{November 29, 2026}


\begin{abstract}
The quantum MacMahon Master Theorem of Garoufalidis, L\^e and
Zeilberger extends the classical identity of MacMahon to right-quantum
matrices, and has since been reproved and extended by several
authors.  We observe that, in the commutative case, the Master Theorem
is a consequence of the matrix-tree theorem of Kirchhoff and Tutte,
and we carry this approach over to the right-quantum algebra.  We prove a quantum matrix-tree
theorem for a generic Manin matrix: the principal minors of its
Laplacian, expanded in the row sums, are again sums over rooted
forests, provided the arcs of each forest are multiplied in a
suitable order.  The order we find depends on the shape of the
forest: it is the order in which Foata and Fuchs encoded forests by
words in their bijective proof of Cayley's formula.  The proof rests
on a noncommutative expansion that removes the path from an arbitrary
vertex to its root.  As a consequence we recover the quantum MacMahon
Master Theorem.
\end{abstract}

\maketitle

\section{Introduction}\label{sec:intro}

The Master Theorem, published by MacMahon in 1915 \cite{MacMahon1960}, is one
of the fundamental results of enumerative combinatorics.  It states that the coefficient of a
monomial in a product of powers of linear forms, whose coefficients
are the rows of a square matrix, is equal to the coefficient of the
same monomial in the inverse of a single determinant.  MacMahon used
it to evaluate sums of products of binomial coefficients, such as
Dixon's identity, and in 1962 Good gave a short proof by Lagrange inversion
\cite{Good1962}.  The algebraic approach goes back to Foata
\cite{Foata1965} and to the 1969 monograph of Cartier and Foata
\cite{Cartier1969Foata}: a product of
entries of the matrix is read as a biword, the row indices forming
the top word and the column indices the bottom word, and the theorem
becomes an identity between formal sums of biwords, valid over a
partially commutative alphabet: the determinant times the sum of the
biwords whose top word is the nondecreasing rearrangement of the
bottom word is equal to one.

In 2006, Garoufalidis, L\^e and Zeilberger \cite{Garoufalidis2006Le} introduced a new kind of
$q$-analogue: the same identity, with a quantum analogue of the
determinant, holds in the \emph{right-quantum algebra}, where the
entries of the matrix no longer commute.  For $q=1$ the right-quantum
matrices are the Manin matrices \cite{Chervov2009Falqui}.  Foata and the author gave a new
proof by determinant calculus, showing that MacMahon's argument
remains valid in the quantum setting \cite{Foata2007Han}, and reduced the
general case to $q=1$ by the ``$1=q$'' principle \cite{Foata2008HanBasis};
Hai and Lorenz derived the theorem from Koszul duality
\cite{Hai2007Lorenz}.  Konvalinka and Pak \cite{Konvalinka2007Pak} gave the first
bijective proof and extended the theorem in
several directions, including a multiparameter analogue.  Their
bijection is related to Foata's first fundamental transformation, a
connection made explicit by Konvalinka \cite{Konvalinka2008}, for total
orders of the alphabet chosen in advance.  In this spirit of \emph{combinatorial linear algebra},
Konvalinka extended further classical identities to right-quantum
matrices: Sylvester's identity \cite{Konvalinka2007}, the inverse matrix
formula, Jacobi's ratio theorem and a generalization of the Master
Theorem \cite{Konvalinka2008}.

\medskip
Another classical determinantal identity is the \emph{matrix-tree
theorem}, going back to Kirchhoff's 1847 work on electrical networks
\cite{Kirchhoff1847}, in the form for directed graphs given by Tutte in
1948 \cite{Tutte1948}: the
principal minors of the Laplacian matrix of a weighted directed graph
are the generating polynomials of its spanning rooted forests.
Among its many extensions, some involve noncommuting quantities: a
large class of combinatorial objects, such as unrooted spanning
forests, can be represented by non-Gaussian Grassmann integrals
\cite{Caracciolo2004}; the edge weights can be quaternions or matrices
\cite{Kassel2020Levy}; and for a generic noncommutative stochastic
matrix, the row vector it fixes, which in the commutative case is given
by the Markov chain tree theorem, has been expressed in terms of codes
\cite{Lavallee2010}.
Somewhat surprisingly, in the commutative case the Master Theorem is a
consequence of the matrix-tree theorem; to the best of our knowledge,
this has not been observed in the literature.  Replace the matrix by
the identity minus the matrix, and multiply by the generating
function of the products of the linear forms.  The matrix-tree
theorem then expresses the determinant times this generating function
as a sum over all rooted forests.  Keep only the balanced terms, in
which every index occurs as often as a row index as a column index.
Every nonempty forest has an arc pointing to a root outside its vertex
set, so only the empty forest survives, and this is exactly the
Master Theorem (Section~\ref{sec:MMT}).

In this paper we prove a matrix-tree theorem in the right-quantum
algebra, that is, for a generic Manin matrix~$\B$, and derive the
quantum Master Theorem from it in the same way: the main theorem
implies the right-quantum Master Theorem at $q=1$, and hence, by the
``$1=q$'' principle, the $q$-version of Garoufalidis, L\^e and
Zeilberger.  A principal minor of
the Laplacian is an alternating sum of principal minors of~$\B$ times
products of row sums.  The row sums commute with each other but not
with the entries of~$\B$; placing them on the right defines the
\emph{forest sums}~$\Delta_J$.  Our main result (Theorem~\ref{thm:main}) states
that each forest sum is again the sum over the rooted forests,
provided the arcs of each forest are multiplied in a suitable order.
The whole difficulty lies in this order: already for two vertices, no
order of the vertices fixed independently of the forest works.  The
correct order, which we call the \emph{regular form}
(Definition~\ref{def:reg}), turns out to be the encoding of forests by
words that Foata and Fuchs introduced in 1970, in their bijective proof
of Cayley's formula \cite{Foata1970Fuchs}.  Roughly, the leaves are
taken in increasing order, and each leaf adds to the regular form the
part of its path towards the root not visited before, written
backwards.  Thus the order
depends on the shape of the forest, and not only on the order of the
alphabet.

The key to the proof is a noncommutative form of the removal of the
path from a vertex to its root (Theorem~\ref{thm:pathintro}), valid
for an arbitrary vertex and proved by absorbing its row sum into the
determinants as an additional column.  This extends the determinantal
argument of \cite{Foata2007Han}, where all the columns are added to
the last one and only the balanced terms of the result are needed; here
the other terms must be computed exactly, and they are the forests in
regular form.  Only the defining relations of
the right-quantum algebra are used; the proof is algebraic, and a
bijective proof remains to be found.  The path expansion also gives a
cycle notation formula for the principal minors of~$\B$, essentially
a formula of Konvalinka \cite{Konvalinka2008} read backwards, and other admissible orderings of
the forests.

Section~\ref{sec:main} gives the precise definitions and statements.
Section~\ref{sec:prelim} contains the commutation rules and the
calculus of determinants whose columns are sums of columns of~$\B$.
Section~\ref{sec:reg} defines the regular form, Section~\ref{sec:path}
proves the path expansion and Section~\ref{sec:proof} the main
theorem.  In Section~\ref{sec:MMT} we derive the quantum MacMahon
Master Theorem from the main theorem.  Section~\ref{sec:cycles}
contains the expansions in terms of cycles, and Section~\ref{sec:cons}
further consequences.

\section{Definitions and main results}\label{sec:main}

\subsection{The matrix-tree theorem}
Let $\setA=\{1,2,\dots,r\}$, and consider the complete directed graph
on~$\setA$.  The weight of the arc $i\to j$ is traditionally denoted
by~$x_{ij}$; here we write it as the biletter $\bw{i}{j}$, with
top~$i$ and bottom~$j$, the notation used in the right-quantum algebra
(Section~\ref{ssec:rq}) and in the rest of the paper.  In this
subsection the $\bw{i}{j}$, $i,j\in\setA$, are commuting
indeterminates.  The \emph{Laplacian} of the graph is the matrix
\[
\mathcal{L}=\Bigl(\delta_{ij}s_i-\bw{i}{j}\Bigr)_{i,j\in\setA},
\qquad s_i=\sum_{a\in\setA}\bw{i}{a},
\]
whose row sums vanish; the loops $\bw{i}{i}$ cancel on the diagonal.  For
$J\subseteq\setA$ we write $\mathcal{L}_J$ for the principal submatrix
with rows and columns in~$J$, and $R=\setA\setminus J$, which will be the set of roots.

The matrix-tree theorem for directed graphs is due to Tutte
\cite{Tutte1948}.  For a vertex $\rho\in\setA$, an \emph{arborescence}
converging to~$\rho$ is a spanning subgraph in which every vertex
other than~$\rho$ has exactly one outgoing arc, $\rho$ has none, and
every vertex has a directed path to~$\rho$.  Tutte's theorem states
that $\det\mathcal{L}_{\setA\setminus\rho}$ is the generating
polynomial of these arborescences, each weighted by the product of
the weights of its arcs.

The \emph{forest version} of the matrix-tree theorem replaces the
single root~$\rho$ by an arbitrary set~$R$ of roots.  A map
$b:J\to\setA$ is \emph{acyclic} if its functional digraph, with arcs
$v\to b(v)$ for $b(v)\in J$, has no cycle.  Equivalently, the arcs
$v\to b(v)$, $v\in J$, form a spanning forest of the complete directed
graph on~$\setA$ in which every vertex of~$J$ has exactly one outgoing
arc and every tree converges to a root in~$R$.

\begin{theorem}[matrix-tree theorem, forest version]\label{thm:forest}
For every $J\subseteq\setA$,
\[
\det\mathcal{L}_J\;=\;\sum_{b:J\to\setA\ \mathrm{acyclic}}\
\prod_{v\in J}\bw{v}{b(v)}.
\]
\end{theorem}

For $R=\{\rho\}$ this is Tutte's theorem.  Conversely, the general
case follows from Tutte's theorem by merging the vertices of~$R$ into
a single root~$\rho$, with arc weights
$\bw{i}{\rho}=\sum_{a\in R}\bw{i}{a}$: this leaves $\mathcal{L}_J$
unchanged, and expanding the products splits each arborescence
converging to~$\rho$ into the forests with roots in~$R$.
Theorem~\ref{thm:forest} is also the principal case of the
\emph{all-minors matrix-tree theorem} of Chaiken \cite{Chaiken1982}, which
expresses an arbitrary minor of~$\mathcal{L}$, obtained by deleting
rows~$W$ and columns~$U$ with $|W|=|U|$, as a signed sum over the
spanning forests in which each tree contains exactly one vertex of~$W$
and one vertex of~$U$; for $W=U=R$ all signs are~$+$.  A
combinatorial proof of Theorem~\ref{thm:forest}, by a sign-reversing
involution on the cycles, is given by Zeilberger \cite{Zeilberger1985};
Kassel and L\'evy revisit it and obtain several variants in a unified
way \cite{Kassel2020Levy}.

Expanding $\det\mathcal{L}_J$ by multilinearity along the diagonal
gives
\begin{equation}\label{eq:forestexp}
\det\mathcal{L}_J\;=\;\sum_{L\subseteq J}(-1)^{|L|}\,
\det\Bigl(\bw{i}{j}\Bigr)_{i,j\in L}\prod_{i\in J\setminus L}s_i\,.
\end{equation}
It is in this form, with the row sums separated from the
determinants, that Theorem~\ref{thm:forest} extends to the right-quantum algebra.

\subsection{The right-quantum algebra}\label{ssec:rq}
A \emph{biletter} is a pair $\bw{x}{a}$ of letters $x,a\in\setA$; $x$
is its \emph{top} and $a$ its \emph{bottom}.  A \emph{biword} is a word
in the biletters, written $\bw{x_1\cdots x_n}{a_1\cdots a_n}
=\bw{x_1}{a_1}\cdots\bw{x_n}{a_n}$.  Let $\calA$ be the free
associative $\mathbb{Z}$-algebra on the biletters.  The \emph{right-quantum
algebra} is $\calR=\calA/\calI$, where $\calI$ is the two-sided
ideal generated by
\begin{align}
&\bw{x\,y}{a\,a}-\bw{y\,x}{a\,a}, \tag{R1}\label{R1}\\[2pt]
&\bw{x\,y}{a\,b}+\bw{x\,y}{b\,a}-\bw{y\,x}{a\,b}-\bw{y\,x}{b\,a},
\tag{R2}\label{R2}
\end{align}
for all $x\neq y$ and $a\neq b$.  We write $E\equiv F$ when
$E-F\in\calI$.  In other words, the generic matrix
$\B=\bigl(\bw{i}{j}\bigr)_{i,j\in\setA}$ is a \emph{Manin matrix}
\cite{Chervov2009Falqui}: entries in a common column commute, and
$[\bw{x}{a},\bw{y}{b}]=[\bw{y}{a},\bw{x}{b}]$.  This is the right-quantum
algebra of \cite{Garoufalidis2006Le} at $q=1$ (Section~\ref{sec:intro}).

For finite sets $T,C\subseteq\setA$ with $|T|=|C|$ and
$C=\{c_1<\dots<c_m\}$, the \emph{column determinant} of $\B_{T,C}$ is
\[
\det(\B_{T,C})=\sum_{\sigma}(-1)^{\inv\sigma}\,
\bw{\sigma_1}{c_1}\bw{\sigma_2}{c_2}\cdots\bw{\sigma_m}{c_m},
\]
where $\sigma=\sigma_1\cdots\sigma_m$ runs over the arrangements of~$T$
and $\inv\sigma$ is the number of inversions; the order of the
factors matters.  Put $\det(\B_K)=\det(\B_{K,K})$ and
$\det(\B_\varnothing)=1$.  The row sums
\[
S_i=\bw{i}{1}+\bw{i}{2}+\cdots+\bw{i}{r}\qquad(i\in\setA)
\]
commute modulo~$\calI$ (Lemma~\ref{lem:gencol}), and we write $S_X$
for their product over a finite set~$X$.  The row sums do not commute
with the biletters, so the determinant of the Laplacian
$\mathrm{diag}(S)-\B$ depends on how the row sums are placed among the
biletters.  Following \eqref{eq:forestexp}, we place them on the
right.  For every finite $J\subseteq\setA$ we define the
\emph{forest sum}~$\Delta_J$ by
\begin{equation}\label{eq:Delta}
\Delta_J\;=\;\sum_{L\subseteq J}(-1)^{|L|}\,\det(\B_L)\,
S_{J\setminus L}\,.
\end{equation}
These are exactly the coefficients of the expansion of the principal
minors of~$\B$ in terms of the row sums, see
Proposition~\ref{prop:mobius}.

\subsection{The main result}
In the commutative image, Theorem~\ref{thm:forest} gives
$\Delta_J=\sum_b\prod_{v\in J}\bw{v}{b(v)}$ over the acyclic maps
$b:J\to\setA$.  In~$\calR$ this formula must be read with care, since
the biletters of a forest no longer commute: the relations
\eqref{R1}--\eqref{R2} allow some reorderings of a word but not others,
and the sum over the forests depends on the order in which the arcs
of each forest are written.  Writing the arcs in increasing order of
their tops already fails for $|J|=2$ (Example~\ref{ex:intro}).
An admissible order must depend on the
shape of the forest, and the one we find is the rearrangement of
Foata and Fuchs \cite{Foata1970Fuchs}.

Let $b:J\to\setA$ be acyclic.  Its \emph{sources} are the letters
of~$J$ that are not values of~$b$, that is, the leaves of the forest.
Every letter of~$J$ is reached from at least one source by
iterating~$b$; assign it to the smallest such source.  The letters
assigned to a source~$u$ form an initial segment
$u,b(u),\dots,b^{s}(u)$ of its orbit, and the \emph{chain} of~$u$ is
the biword
$\bw{b^{s}(u)}{b^{s+1}(u)}\cdots\bw{b(u)}{b^2(u)}\bw{u}{b(u)}$,
the segment written backwards, each letter carrying its image.  The
\emph{regular form} $\Reg(b)$ is the product of the chains of the
sources, in increasing order of the sources
(Definition~\ref{def:reg} and Example~\ref{ex:tree}).  In $\Reg(b)$
every letter appears in the top row before its occurrences in the
bottom row.  The bottom row of $\Reg(b)$ is the word that Foata and
Fuchs associate with the forest (Remark~\ref{rem:FF}).

Our main theorem, proved in Section~\ref{sec:proof}, is stated next.

\begin{theorem}[right-quantum matrix-tree theorem]\label{thm:main}
For every finite $J\subseteq\setA$,
\[
\Delta_J\;\equiv\;\sum_{b:J\to\setA\ \mathrm{acyclic}}\Reg(b).
\]
Equivalently, for every finite $K\subseteq\setA$,
\[
\det(\B_K)\;\equiv\;\sum_{J\subseteq K}(-1)^{|J|}\,
\Bigl(\sum_{b:J\to\setA\ \mathrm{acyclic}}\Reg(b)\Bigr)\,
S_{K\setminus J}\,.
\]
\end{theorem}

In the commutative image, and more generally in the Cartier--Foata quotient
\cite{Cartier1969Foata,Foata1965}, where biletters with distinct tops commute,
$\Reg(b)$ may be replaced by the product of the biletters of~$b$ in
any order, and Theorem~\ref{thm:main} reduces to
Theorem~\ref{thm:forest}.

\begin{example}\label{ex:intro}
Let $r=3$ and $J=\{1,2\}$.  The acyclic maps on~$J$ are
$b=(b(1),b(2))=(2,3)$, $(3,1)$, $(3,3)$, with regular forms
$\bw{21}{32}$, $\bw{12}{31}$, $\bw{12}{33}$, so that
\[
\Delta_{\{1,2\}}\;\equiv\;\bw{21}{32}+\bw{12}{31}+\bw{12}{33}\,.
\]
The three maps have distinct bottom contents and $\calI$ is graded by
the bottom content, so each term must be written correctly by itself.
The forest $1\to2\to3$ appears as $\bw{21}{32}$, with the letter~$2$
first, and the other arrangement $\bw{12}{23}$ is not congruent to
it: in the basis of~$\calR$ formed by the biwords without double
descent \cite{Foata2008HanBasis}, $\bw{12}{23}$ is a basis element while
$\bw{21}{32}\equiv\bw{12}{32}+\bw{12}{23}-\bw{21}{23}$.  The forest
$2\to1\to3$, on the contrary, appears as $\bw{12}{31}$, with the
letter~$1$ first, and the other arrangement $\bw{21}{13}$ is a basis
element different from $\bw{12}{31}$.  Hence no order of the vertices
chosen independently of the forest is admissible, already for
$|J|=2$.
\end{example}

\subsection{The path expansion}
A natural way to decompose a forest is to remove the path from a
given vertex to its root.  The key to
Theorem~\ref{thm:main} is a noncommutative form of this
decomposition.  Let $u\in J$.  An \emph{escaping path} from~$u$
in~$J$ is a pair $(P,\beta)$, where $P=(u,v_1,\dots,v_s)$ is a sequence
of distinct letters of~$J$ ($s\ge0$) and $\beta\in\setA\setminus J$;
its \emph{support} is $V=\{u,v_1,\dots,v_s\}$ and its \emph{word} is
\begin{equation}\label{eq:wP}
w(P,\beta)\;=\;\bw{v_s}{\beta}\bw{v_{s-1}}{v_s}\cdots
\bw{v_1}{v_2}\bw{u}{v_1}
\qquad(\text{for } s=0:\ w(P,\beta)=\bw{u}{\beta}).
\end{equation}

\begin{theorem}[path expansion]\label{thm:pathintro}
Let $J\subseteq\setA$ be finite and let $u$ be any letter of~$J$.
Then
\[
\Delta_J\;\equiv\;\sum_{(P,\beta)}w(P,\beta)\,\Delta_{J\setminus V}\,,
\]
where $(P,\beta)$ runs over the escaping paths from~$u$ in~$J$, of
support~$V$.
\end{theorem}

\begin{example}\label{ex:path}
Take $r=3$ and $J=\{1,2\}$ as in Example~\ref{ex:intro}, so that the
only possible end is $\beta=3$.  For a single letter,
$\Delta_{\{a\}}=S_a-\bw{a}{a}$ is the sum of the biletters $\bw{a}{c}$
with $c\neq a$.  For $u=1$ the escaping paths are $P=(1)$, with word
$\bw{1}{3}$ and remaining set $\{2\}$, and $P=(1,2)$, with word
$\bw{21}{32}$ and nothing left.  Hence
\[
\Delta_{\{1,2\}}\;\equiv\;\bw{1}{3}\,\Delta_{\{2\}}+\bw{21}{32}
\;=\;\bw{12}{31}+\bw{12}{33}+\bw{21}{32}\,,
\]
which are exactly the regular forms of Example~\ref{ex:intro}.  For
$u=2$ the escaping paths are $P=(2)$ and $P=(2,1)$, and
\[
\Delta_{\{1,2\}}\;\equiv\;\bw{2}{3}\,\Delta_{\{1\}}+\bw{12}{31}
\;=\;\bw{21}{32}+\bw{21}{33}+\bw{12}{31}\,.
\]
The two expansions agree, since $\bw{21}{33}\equiv\bw{12}{33}$
by~\eqref{R1}.  In each case the path from~$u$ is written in front,
read from its end~$\beta$ back to~$u$, and the forests of the remaining
letters follow.
\end{example}

The path is written in front, in the order \eqref{eq:wP}, for an
arbitrary~$u$, and no condition is imposed on the remaining factor.
The proof (Section~\ref{sec:path}) absorbs the row sum~$S_u$, which
stands on the right in~$\Delta_J$, into the determinants as an
additional column.  Theorem~\ref{thm:main} does not follow from
Theorem~\ref{thm:pathintro} directly, since the first chain of
$\Reg(b)$ starts at the smallest source, which depends on~$b$.  The
two are connected through a word-for-word factorization of the
regular forms at their first chain (Proposition~\ref{prop:firstchain})
and a decomposition according to the bottom letters, the slices of
Definition~\ref{def:slice} (Lemma~\ref{lem:slices} and
Section~\ref{sec:proof}).

\subsection{Consequences}
The substitution $\B\mapsto\I-\B$ turns Theorem~\ref{thm:main} into
the quantum MacMahon Master Theorem of~\cite{Garoufalidis2006Le}
(Corollary~\ref{cor:MMT}).  Further, the path expansion contains a
cycle expansion of the determinant at
an arbitrary letter (Corollary~\ref{cor:cyclepeel}) and, by iteration,
a cycle notation formula for $\det(\B_M)$
(Corollary~\ref{cor:cycleform}), which is essentially a formula of
Konvalinka \cite{Konvalinka2008} read backwards (Remark~\ref{rem:Konv}).
Theorem~\ref{thm:main} and this formula are the two extreme cases of a
single identity for the sums
$\sum_{L\subseteq J}(-1)^{|L|}\det(\B_{L\cup M})S_{J\setminus L}$
(Theorem~\ref{thm:master}).  The regular form is not the only
admissible ordering: any rule that peels the forest along the orbit of
a letter chosen in advance also works (Corollary~\ref{cor:selector}).  Keeping only the bottom
row of the biwords turns Theorem~\ref{thm:main} into the Foata--Fuchs
encoding of forests, and hence into Cayley's formula
(Corollary~\ref{cor:cayley}).

\section{Basic properties of the right-quantum algebra}\label{sec:prelim}

\subsection{Contents}
The \emph{top content} and the \emph{bottom content} of a biword are
the multisets of its top and bottom letters.  Each generator
\eqref{R1}--\eqref{R2} of~$\calI$ is a combination of biwords with the
same top content and the same bottom content.  Hence, if $\pi$ is the
linear projection of~$\calA$ that keeps the biwords whose pair (top
content, bottom content) belongs to a given set and kills the others,
then $\pi(\calI)\subseteq\calI$; in other words, congruences may be
projected:
\begin{equation}\label{eq:piI}
E\equiv F\quad\Longrightarrow\quad\pi(E)\equiv\pi(F).
\end{equation}

\subsection{Commutation rules}
For letters $x\neq y$ and $a,b$ put
$\gamma_{x,y}(a,b)=\bw{x}{a}\bw{y}{b}-\bw{y}{a}\bw{x}{b}$.  Then
\eqref{R1} is $\gamma_{x,y}(a,a)$ and \eqref{R2} is
$\gamma_{x,y}(a,b)+\gamma_{x,y}(b,a)$, so that
\begin{equation}\label{eq:gamma}
\gamma_{x,y}(a,b)\;\equiv\;-\gamma_{x,y}(b,a)
\quad\text{and}\quad \gamma_{x,y}(a,a)\;\equiv\;0
\qquad(x\neq y).
\end{equation}

For a letter $x$ and a set $X\subseteq\setA$ put
$\bw{x}{X}=\sum_{a\in X}\bw{x}{a}$, so that $S_x=\bw{x}{\setA}$.

\begin{lemma}\label{lem:gencol}
Let $x\neq y$ be letters and $X,Y\subseteq\setA$.  Then
\begin{align}
\bw{x}{X}\bw{y}{Y}-\bw{y}{X}\bw{x}{Y}
&\;\equiv\;-\Bigl(\bw{x}{Y}\bw{y}{X}-\bw{y}{Y}\bw{x}{X}\Bigr),
\label{eq:gexch}\\
\bw{x}{X}\bw{y}{X}&\;\equiv\;\bw{y}{X}\bw{x}{X}.
\label{eq:gcomm}
\end{align}
In particular $S_xS_y\equiv S_yS_x$ for all letters $x,y$.
\end{lemma}

\begin{proof}
With $\gamma=\gamma_{x,y}$, the two sides of \eqref{eq:gexch} are
$\sum_{a\in X,\,b\in Y}\gamma(a,b)$ and
$-\sum_{a\in X,\,b\in Y}\gamma(b,a)$, which are congruent term by term
by \eqref{eq:gamma}.  The difference of the two sides of
\eqref{eq:gcomm} is
\[
\sum_{a\in X}\gamma(a,a)+\sum_{a,b\in X,\ a<b}
\bigl(\gamma(a,b)+\gamma(b,a)\bigr),
\]
and each term is $\equiv0$ by \eqref{eq:gamma}.  Taking $X=Y=\setA$
in \eqref{eq:gcomm} gives $S_xS_y\equiv S_yS_x$ for $x\neq y$, and the
case $x=y$ is trivial.
\end{proof}

\subsection{Determinants in the right-quantum algebra}
For a matrix with commuting entries, the determinant has many simple
properties: it is additive in each column, it changes sign when two
columns are exchanged, it vanishes when two columns are equal, and it
can be expanded along a column.  We show that these properties remain
true for column determinants in the right-quantum algebra, provided
they are handled with care: the exchange of columns holds only modulo
the ideal~$\calI$ and for columns of the form~$[X]$ defined below, and
in an expansion along the last (resp.\ first) column, the entries of
that column stand on the right (resp.\ left).

Let $T=\{t_1<\dots<t_m\}$ be a finite set of letters.  A \emph{column}
is a family $\gamma=(\gamma(t))_t$ of expressions indexed by letters,
and for columns $\gamma_1,\dots,\gamma_m$ we put
\begin{equation}\label{eq:gendet}
\det(T;\gamma_1,\dots,\gamma_m)\;=\;
\sum_{\sigma}(-1)^{\inv\sigma}\,
\gamma_1(\sigma_1)\gamma_2(\sigma_2)\cdots\gamma_m(\sigma_m),
\end{equation}
the sum being over the arrangements $\sigma=\sigma_1\cdots\sigma_m$
of~$T$ ($\det(\varnothing;)=1$).  For $X\subseteq\setA$ let $[X]$ be
the column $t\mapsto\bw{t}{X}$, and $[a]=[\{a\}]$.  Thus
$\det(\B_{T,C})=\det(T;[c_1],\dots,[c_m])$ for $C=\{c_1<\dots<c_m\}$,
and $[\setA]$ is the column of the row sums.  For $t\in T$ let
$\rk_T(t)=\#\{t'\in T:\ t'\le t\}$ be the rank of~$t$ in~$T$.

\begin{lemma}\label{lem:gendet}
\textup{(i)} The determinant \eqref{eq:gendet} is additive in each
column; in particular
$\det(T;\dots,[X\cup Y],\dots)
=\det(T;\dots,[X],\dots)+\det(T;\dots,[Y],\dots)$ for disjoint
$X,Y\subseteq\setA$.

\textup{(ii)} If all the columns are of the form $[X]$, exchanging two
adjacent columns changes the sign modulo~$\calI$.

\textup{(iii)} If all the columns are of the form $[X]$ and two of
them are equal, the determinant is $\equiv0$.

\textup{(iv)} For arbitrary columns,
\[
\det(T;\gamma_1,\dots,\gamma_m)
=\sum_{t\in T}(-1)^{m-\rk_T(t)}\,
\det(T\setminus t;\gamma_1,\dots,\gamma_{m-1})\;\gamma_m(t).
\]

\textup{(v)} For arbitrary columns,
\[
\det(T;\gamma_1,\dots,\gamma_m)
=\sum_{t\in T}(-1)^{\rk_T(t)-1}\,\gamma_1(t)\;
\det(T\setminus t;\gamma_2,\dots,\gamma_m).
\]
\end{lemma}

\begin{proof}
(i) is clear.  (ii) Let the two columns $[X],[Y]$ stand in positions
$i,i+1$.  Pair each arrangement~$\sigma$ with
$\sigma_i<\sigma_{i+1}$ with the arrangement obtained by exchanging
these two letters.  With $x=\sigma_i$, $y=\sigma_{i+1}$, and $p$, $q$
the products of the factors before and after the two positions, the
pair contributes
$(-1)^{\inv\sigma}p\bigl(\bw{x}{X}\bw{y}{Y}-\bw{y}{X}\bw{x}{Y}\bigr)q$
to the determinant and
$(-1)^{\inv\sigma}p\bigl(\bw{x}{Y}\bw{y}{X}-\bw{y}{Y}\bw{x}{X}\bigr)q$
to the determinant with the two columns exchanged.  The sum is
$\equiv0$ by \eqref{eq:gexch}.  (iii) For two equal adjacent columns
the same pairing gives contributions
$\pm p\bigl(\bw{x}{X}\bw{y}{X}-\bw{y}{X}\bw{x}{X}\bigr)q\equiv0$ by
\eqref{eq:gcomm}; the general case follows by~(ii).  (iv) and (v):
group the arrangements according to their last, respectively first,
letter~$t$, which is involved in $m-\rk_T(t)$, respectively
$\rk_T(t)-1$, inversions.
\end{proof}

Properties (i) and (iv) are \cite[Properties~4 and~5]{Foata2007Han},
property (ii) for the columns of~$\B$ is \cite[Property~6]{Foata2007Han},
and the commutation of the row sums in Lemma~\ref{lem:gencol} is
\cite[Lemma~9]{Foata2007Han}; (ii) and (iii) reflect the fact
that~$\B$ is a Manin matrix.

\subsection{Reversal of biwords}\label{ssec:rev}
Let $\omega$ be the linear map of~$\calA$ that reverses every biword,
$\omega\bigl(\bw{x_1}{a_1}\cdots\bw{x_n}{a_n}\bigr)
=\bw{x_n}{a_n}\cdots\bw{x_1}{a_1}$; it is an anti-automorphism.

\begin{lemma}\label{lem:rev}
$\omega(\calI)=\calI$, $\omega(S_i)=S_i$, and
$\omega(\det\B_M)\equiv\det(\B_M)$ for every finite
$M\subseteq\setA$.
\end{lemma}

\begin{proof}
$\omega$ maps the generators \eqref{R1} and \eqref{R2} to their
opposites, and fixes each biletter.  Let $M=\{c_1<\dots<c_m\}$.  The
reversal of the term of $\sigma=\sigma_1\cdots\sigma_m$ in
$\det(\B_M)$ is $\bw{\sigma_m}{c_m}\cdots\bw{\sigma_1}{c_1}$, the term
of the arrangement $\sigma_m\cdots\sigma_1$, which has
$\binom m2-\inv\sigma$ inversions, in $\det(M;[c_m],\dots,[c_1])$.
Hence $\omega(\det\B_M)=(-1)^{\binom m2}\det(M;[c_m],\dots,[c_1])$,
and reversing the order of the $m$ columns by $\binom m2$ adjacent
exchanges gives $\det(\B_M)$ back, by
Lemma~\ref{lem:gendet}\,(ii).
\end{proof}

\section{Acyclic maps and the regular form}\label{sec:reg}

Let $J\subseteq\setA$ be finite.  A \emph{map on~$J$} is a function
$b:J\to\setA$; its biletters are $\bw{v}{b(v)}$, $v\in J$, and any
biword formed by these $|J|$ biletters is an \emph{arrangement}
of~$b$.  The map~$b$ is \emph{acyclic} if its functional digraph
(arcs $v\to b(v)$ for $b(v)\in J$) has no cycle, that is, if starting
from any $v\in J$ and applying~$b$ repeatedly, one eventually reaches a
letter outside~$J$, a root.  Let $b$ be acyclic.
Its set of \emph{sources} is $\Src(b)=J\setminus b(J)$; it is
nonempty if $J\neq\varnothing$.  For $v\in J$ let
\[
\kappa(v)=\min\{u\in\Src(b):\ v\in\{u,b(u),b^2(u),\dots\}\}.
\]
Along a forward orbit~$\kappa$ is non-increasing, and a source~$u$
satisfies $\kappa(u)=u$; hence $\kappa^{-1}(u)$ is an initial segment
$u,b(u),\dots,b^{s}(u)$ of the orbit of~$u$.

\begin{definition}\label{def:reg}
Let $b$ be acyclic on~$J$ with sources $u_1<\dots<u_t$, and let
$\kappa^{-1}(u_i)=\{u_i,b(u_i),\dots,b^{s_i}(u_i)\}$.  Then the
\emph{chain} of~$u_i$ is
\[
\ch(u_i)\;=\;
\bw{b^{s_i}(u_i)}{b^{s_i+1}(u_i)}\cdots
\bw{b(u_i)}{b^{2}(u_i)}\,\bw{u_i}{b(u_i)},
\]
and the \emph{regular form} of~$b$ is
\begin{equation}\label{eq:Reg}
\Reg(b)\;=\;\ch(u_1)\,\ch(u_2)\cdots\ch(u_t).
\end{equation}
We put
\begin{equation}\label{eq:FJ}
\Gamma_J\;=\;\sum_{b:J\to\setA\ \mathrm{acyclic}}\Reg(b),
\qquad \Gamma_\varnothing=1.
\end{equation}
\end{definition}

\begin{example}\label{ex:tree}
Let $J=\{1,3,5,6,7,9,10,12,13,14\}$ and let $b$ be the acyclic map
\[
\begin{array}{c|cccccccccc}
v&1&3&5&6&7&9&10&12&13&14\\ \hline
b(v)&14&11&14&12&1&11&1&8&6&12
\end{array}
\]
with values in $\setA=\{1,\dots,14\}$.  Its forest is drawn below,
each arrow pointing from a letter to its image; the two roots $11$ and
$8$, which lie outside~$J$, are boxed.
\begin{center}
\begin{tikzpicture}[>=Stealth,scale=0.95,every node/.style={font=\small}]
\node (n11) at (0,0) [draw,rectangle] {$11$};
\node (n3) at (-0.9,1.1) {$3$};
\node (n9) at (0.9,1.1) {$9$};
\draw[->] (n3) -- (n11);
\draw[->] (n9) -- (n11);
\node (n8) at (4,0) [draw,rectangle] {$8$};
\node (n12) at (4,1.1) {$12$};
\draw[->] (n12) -- (n8);
\node (n14) at (3,2.2) {$14$};
\node (n6) at (5.3,2.2) {$6$};
\draw[->] (n14) -- (n12);
\draw[->] (n6) -- (n12);
\node (n1) at (2.2,3.3) {$1$};
\node (n5) at (3.8,3.3) {$5$};
\draw[->] (n1) -- (n14);
\draw[->] (n5) -- (n14);
\node (n13) at (5.3,3.3) {$13$};
\draw[->] (n13) -- (n6);
\node (n7) at (1.5,4.4) {$7$};
\node (n10) at (2.9,4.4) {$10$};
\draw[->] (n7) -- (n1);
\draw[->] (n10) -- (n1);
\end{tikzpicture}
\end{center}
The sources (leaves) are $3<5<7<9<10<13$, and their chains are formed
in this order.  The source~$5$ is the smallest one in the tree of the
root~$8$, so its chain takes the whole path $5\to14\to12\to8$; written
backwards, it is $\bwc{12/8,14/12,5/14}$.  The path of the source~$13$
is $13\to6\to12\to8$, but $12$ has already been taken by the chain
of~$5$; so the chain of~$13$ keeps only the initial segment
$13\to6\to12$ of this path, and is $\bwc{6/12,13/6}$.  In the same way,
the chain of~$7$ stops at~$14$ and the chain of~$10$ stops at~$1$.  In
terms of~$\kappa$, this means $\kappa(14)=\kappa(12)=5$, $\kappa(1)=7$
and $\kappa(6)=13$.  Hence
\[
\Reg(b)=\underbrace{\bw{3}{11}}_{\ch(3)}\;
\underbrace{\bwc{12/8,14/12,5/14}}_{\ch(5)}\;
\underbrace{\bwc{1/14,7/1}}_{\ch(7)}\;
\underbrace{\bw{9}{11}}_{\ch(9)}\;
\underbrace{\bw{10}{1}}_{\ch(10)}\;
\underbrace{\bwc{6/12,13/6}}_{\ch(13)}\,.
\]
\end{example}

\begin{remark}\label{rem:FF}
The regular form is the encoding of forests of Foata and Fuchs
\cite{Foata1970Fuchs}.  Their bijection~$\phi$ maps every map~$f$ of a
finite set into itself to a rearrangement $\phi(f)$ of the word of its
values \cite[Theorem~5.9]{Foata1970Fuchs}: $\phi(f)=\sigma_0\sigma_1\cdots
\sigma_p$, where $\sigma_0$ encodes the cycles of~$f$, $z_1<\dots<z_p$
are the elements that are not values of~$f$, and
$\sigma_k=f^{m_k}(z_k)\cdots f^2(z_k)f(z_k)$, with $m_k$ the least
$m>0$ such that $f^m(z_k)$ is a letter of
$\sigma_0\sigma_1\cdots\sigma_{k-1}$ \cite[(5.2)--(5.5)]{Foata1970Fuchs}.
Extend an acyclic map~$b$ on~$J$ to the map~$f$ of~$\setA$ into itself
that fixes every letter of $\setA\setminus J$.  Then the $z_k$ are the
sources of~$b$, $\sigma_0$ is the increasing word of the letters of
$\setA\setminus J$, and $\sigma_k$ is the bottom row of the chain
of~$z_k$; hence $\phi(f)$ is $\sigma_0$ followed by the bottom row of
$\Reg(b)$.  In Example~\ref{ex:tree}, $\sigma_0=2\,4\,8\,11$ and the
bottom row of $\Reg(b)$ is $11\;8\,12\,14\;14\,1\;11\;1\;12\,6$.  Since
$\phi$ is a bijection, the bottom row of $\Reg(b)$ determines~$b$;
Foata and Fuchs use this to count forests and to recover Cayley's
formula \cite[Section~6]{Foata1970Fuchs}.  For a map with cycles they
first write the cycles, as in the Cartier--Foata theory of
permutations; Theorem~\ref{thm:main} concerns the acyclic part.
\end{remark}

\begin{definition}[slice]\label{def:slice}
Let $G\subseteq\setA$ be finite, $u\in\setA$, and let $E\in\calA$ be
such that all the biwords of~$E$ have top content~$G$.  The
\emph{slice} of~$E$ at~$u$, denoted $E^{\langle u\rangle}$, is the
projection of~$E$ on the biwords whose bottom letters avoid~$u$ and
cover $\{j\in G:\ j<u\}$.
\end{definition}

The condition depends only on the bottom content, so
by~\eqref{eq:piI}, $E\equiv F$ implies
$E^{\langle u\rangle}\equiv F^{\langle u\rangle}$.

\begin{example}\label{ex:projE}
Let $r=8$, $G=\{1,3,4,5,7\}$ and
\[
E=\bw{31547}{34776}+\bw{17453}{15332}+\bw{31547}{18866}
-\bw{17453}{34521}+\bw{31547}{54138}.
\]
Every biword of~$E$ has top content~$G$.  For $u=2$ the bottom row
must avoid~$2$ and contain~$1$; for $u=3$, a letter of~$G$, it must
avoid~$3$ and contain~$1$; for $u=6$ it must avoid~$6$ and contain
$1,3,4,5$.  Hence
\begin{gather*}
E^{\langle 2\rangle}=\bw{31547}{18866}+\bw{31547}{54138},\qquad
E^{\langle 3\rangle}=\bw{31547}{18866},\\
E^{\langle 6\rangle}=-\bw{17453}{34521}+\bw{31547}{54138}.
\end{gather*}
\end{example}

Every biword of $\Gamma_G$ has top content~$G$.  Since the bottom
letters of $\Reg(b')$ are the values of~$b'$, for $u\notin G$ the slice
$(\Gamma_G)^{\langle u\rangle}$ is the sum of $\Reg(b')$ over the acyclic
maps~$b'$ on~$G$ such that $u\notin b'(G)$ and
$\{j\in G:\ j<u\}\subseteq b'(G)$.  Escaping paths and their words are
as in \eqref{eq:wP}.

\begin{proposition}[first-chain factorization]\label{prop:firstchain}
For every nonempty finite $K\subseteq\setA$, in the free
algebra~$\calA$,
\[
\Gamma_K\;=\;\sum_{u\in K}\ \sum_{(P,\beta)}w(P,\beta)\;
(\Gamma_{K\setminus V})^{\langle u\rangle}\,,
\]
where $(P,\beta)$ runs over the escaping paths from~$u$ in~$K$, of
support~$V$.
\end{proposition}

\begin{proof}
We give a bijection between the acyclic maps~$b$ on~$K$ and the
quadruples $(u,P,\beta,b')$ indexing the right-hand side, where
$(P,\beta)$ is an escaping path from~$u$ in~$K$, of support~$V$, and
$b'$ is an acyclic map on $K\setminus V$ with
$u\notin b'(K\setminus V)\supseteq\{j\in K\setminus V:\ j<u\}$, such
that $\Reg(b)=w(P,\beta)\Reg(b')$.  In words, the bijection cuts
$\Reg(b)=\ch(u_1)\ch(u_2)\cdots\ch(u_t)$ after its first chain:
$\ch(u_1)=w(P,\beta)$ and $\ch(u_2)\cdots\ch(u_t)=\Reg(b')$.  We have
to check that this cut is well defined, that the resulting~$b'$
satisfies the two conditions on its values, and that every quadruple
satisfying them arises exactly once.

Let $b$ be acyclic on~$K$ and $u=\min\Src(b)$.  Since $\kappa\ge u$
and $\kappa$ is non-increasing along orbits, $\kappa=u$ on the whole
orbit of~$u$; hence $\ch(u)=w(P,\beta)$ with
$P=(u,b(u),\dots,b^s(u))$ and $\beta=b^{s+1}(u)\notin K$.  Let $V$ be
the support of~$P$ and $b'=b|_{K\setminus V}$.  The letters of~$V$ are mapped into
$V\cup\{\beta\}$, so a path of~$b$ that enters~$V$ never comes back
to $K\setminus V$.  Hence the letters of $K\setminus V$ have the same
preimages, sources and paths from their sources for $b$ and~$b'$: the
chains of~$b'$ are the other chains of~$b$, and
$\Reg(b)=w(P,\beta)\Reg(b')$.  The conditions hold since $u$ is a
source of~$b$, and every $j<u$ in $K\setminus V$, not being a source
of~$b$, is a value of~$b$, hence of~$b'$.

Conversely, given a quadruple, let $b$ extend~$b'$ by following~$P$
to~$\beta$.  It is acyclic, and by the two conditions its sources
are~$u$ and the sources of~$b'$, all larger than~$u$.  Hence
$u=\min\Src(b)$, and the first construction applied to~$b$ gives back
the quadruple.
\end{proof}

\begin{lemma}\label{lem:slices}
Let $K\subseteq\setA$ be finite and nonempty, and let $\pi_*$ be the
projection on the biwords whose bottom letters cover~$K$.
\textup{(a)} For every $E\in\calA$ all of whose biwords have top
content~$K$,
$E=\pi_*(E)+\sum_{v\in K}E^{\langle v\rangle}$.
\textup{(b)} Let $(P,\beta)$ be an escaping path from~$u$ in~$K$, of
support~$V$, and let $\xi\in\calA$ be such that all the biwords
of~$\xi$ have top content $K\setminus V$.  Then
$\bigl(w(P,\beta)\,\xi\bigr)^{\langle u\rangle}
=w(P,\beta)\,\xi^{\langle u\rangle}$.
\end{lemma}

\begin{proof}
(a) A set of bottom letters that covers~$K$ satisfies the condition
of the slice at~$v$ for no $v\in K$.  One that does not cover~$K$
satisfies it for exactly one $v\in K$, namely the smallest letter
of~$K$ it misses.

(b) The bottom letters of $w(P,\beta)$, $P=(u,v_1,\dots,v_s)$, are
$v_1,\dots,v_s,\beta$; they differ from~$u$, and
$\{j\in K:\ j<u\}\setminus\{v_1,\dots,v_s\}=\{j\in K\setminus V:\ j<u\}$.
Hence a biword $w(P,\beta)\,\xi'$ with $\xi'$ of top content
$K\setminus V$ is kept by the slice at~$u$ if and only if the bottom
letters of~$\xi'$ avoid~$u$ and cover $\{j\in K\setminus V:\ j<u\}$.
\end{proof}

\begin{example}\label{ex:proj}
Let $r=4$ and $K=\{1,2,3\}$, and consider the terms of
Proposition~\ref{prop:firstchain} with $u=2$, $P=(2)$ and $\beta=4$,
so that $K\setminus V=\{1,3\}$.  The slice at~$2$ of
$\Gamma_{\{1,3\}}$ keeps the biwords whose bottom row contains~$1$
but not~$2$.  By~\eqref{eq:FJ}, summing $\Reg(b')$ over the eight
acyclic maps $b':\{1,3\}\to\setA$, we get
\[
\Gamma_{\{1,3\}}=\bw{13}{21}+\bw{13}{22}+\bw{13}{24}+\bw{31}{23}
+\bw{31}{43}+\bw{13}{41}+\bw{13}{42}+\bw{13}{44},
\]
so that $(\Gamma_{\{1,3\}})^{\langle 2\rangle}=\bw{13}{41}$, and the
only such term is
\[
\bw{2}{4}\,(\Gamma_{\{1,3\}})^{\langle 2\rangle}=\bw{213}{441}=\Reg(b),
\qquad (b(1),b(2),b(3))=(4,4,1),
\]
the smallest source of~$b$ being~$2$.  If instead the same map~$b$ is
split at the letter~$1$, with $P=(1)$ and $\beta=4$, we get the word
$\bw{1}{4}\Reg(b|_{\{2,3\}})=\bw{123}{441}$, which is not $\Reg(b)$.
Its bottom letters miss~$2$, the smallest letter of~$K$ they miss, so
it belongs to the slice at~$2$ (Lemma~\ref{lem:slices}\,(a)); this is
how the proof of Theorem~\ref{thm:main} in Section~\ref{sec:proof}
deals with such words.
\end{example}

\section{The path expansion}\label{sec:path}

Throughout this section $J\subseteq\setA$ is finite, $u\in J$, and
$J'=J\setminus\{u\}$.  For $L=\{l_1<\dots<l_k\}\subseteq J'$ let
$\nu_L=\#\{l\in L:\ l>u\}$, and for $X\subseteq\setA$ put
\begin{equation}\label{eq:NL}
N_L(X)\;=\;(-1)^{\nu_L}\,
\det\bigl(L\cup u;\,[l_1],\dots,[l_k],[X]\bigr),
\end{equation}
and $N_L(a)=N_L(\{a\})$.  Up to sign, this is the determinant with
rows $L\cup u$ formed by the columns of~$\B$ indexed by~$L$ and the
column~$[X]$.

\begin{lemma}\label{lem:Nprops}
Let $L\subseteq J'$.
\textup{(a)} $N_L(X\cup Y)=N_L(X)+N_L(Y)$ for disjoint $X,Y$.
\textup{(b)} $N_L(a)\equiv0$ for $a\in L$.
\textup{(c)} $N_L(u)\equiv\det(\B_{L\cup u})$.
\textup{(d)} For every $x\in L$,
\[
\det\bigl((L\cup u)\setminus x;\,[l_1],\dots,[l_k]\bigr)
\;\equiv\;(-1)^{\nu_{L\setminus x}+\#\{l\in L:\,l>x\}}\,
N_{L\setminus x}(x).
\]
\end{lemma}

\begin{proof}
(a) and (b) are Lemma~\ref{lem:gendet}\,(i) and~(iii).  (c) Moving
the column $[u]$ to the left across the $\nu_L$ columns $[l]$ with
$l>u$ costs $(-1)^{\nu_L}$ (Lemma~\ref{lem:gendet}\,(ii)) and
produces $\det(\B_{L\cup u})$.  (d) Moving the column $[x]$ to the
last position costs $(-1)^{\#\{l\in L:\,l>x\}}$ and produces
\[
\det\bigl((L\setminus x)\cup u;\dots,\widehat{[x]},\dots,[x]\bigr)
=(-1)^{\nu_{L\setminus x}}N_{L\setminus x}(x).\qedhere
\]
\end{proof}

\begin{lemma}\label{lem:Nexpan}
For every $L\subseteq J'$ and every $X\subseteq\setA$,
\begin{align}
N_L(X)&\;\equiv\;\det(\B_L)\,\bw{u}{X}
 \;-\;\sum_{x\in L}N_{L\setminus x}(x)\,\bw{x}{X},
\tag{E$_r$}\label{eq:Er}\\
N_L(X)&\;\equiv\;\bw{u}{X}\,\det(\B_L)
 \;-\;\sum_{x\in L}\bw{x}{X}\,N_{L\setminus x}(x).
\tag{E$_l$}\label{eq:El}
\end{align}
\end{lemma}

\begin{proof}
Let $k=|L|$ and $T=L\cup u$.  Then $\rk_T(u)=1+k-\nu_L$ and
$\rk_T(x)=k-\#\{l\in L:\ l>x\}+[u<x]$ for $x\in L$, where $[u<x]$ is
the Iverson bracket \cite[Section~2.1]{Graham1994KP}, equal to $1$
if $u<x$ and $0$ otherwise; moreover
$\nu_L-\nu_{L\setminus x}=[u<x]$.

\eqref{eq:Er}: expand the determinant in \eqref{eq:NL} by
Lemma~\ref{lem:gendet}\,(iv).  The term $t=u$ is
$(-1)^{\nu_L}\det(\B_L)\bw{u}{X}$.  The term $t=x\in L$ is
\begin{align*}
&(-1)^{1+\#\{l>x\}-[u<x]}\det(T\setminus x;[l_1],\dots,[l_k])\bw{x}{X}\\
&\qquad\equiv\;
(-1)^{1-[u<x]+\nu_{L\setminus x}}N_{L\setminus x}(x)\bw{x}{X}
\;=\;-(-1)^{\nu_L}N_{L\setminus x}(x)\bw{x}{X}
\end{align*}
by Lemma~\ref{lem:Nprops}\,(d).

\eqref{eq:El}: by Lemma~\ref{lem:gendet}\,(ii), moving the column
$[X]$ to the first position costs $(-1)^k$; then expand by
Lemma~\ref{lem:gendet}\,(v).  The term $t=u$ carries
the sign $(-1)^{k+\rk_T(u)-1}=(-1)^{\nu_L}$.  The term $t=x$ carries
the sign
\[
(-1)^{k+\rk_T(x)-1}=(-1)^{1+\#\{l>x\}+[u<x]}
\]
and, by Lemma~\ref{lem:Nprops}\,(d), is congruent to
\[
(-1)^{1+[u<x]+\nu_{L\setminus x}}\bw{x}{X}N_{L\setminus x}(x)
=-(-1)^{\nu_L}\bw{x}{X}N_{L\setminus x}(x).
\]
In both cases every term of the expansion is $(-1)^{\nu_L}$ times the
corresponding term of the right-hand side.  Since $N_L(X)$ is
$(-1)^{\nu_L}$ times the determinant, multiplying by $(-1)^{\nu_L}$
gives \eqref{eq:Er} and \eqref{eq:El}.
\end{proof}

Lemmas~\ref{lem:Nprops} and~\ref{lem:Nexpan} extend the determinant
calculus of \cite{Foata2007Han}.  There, for the matrix $\I-\B$ on the
whole alphabet, all the columns are added to the last one; the
determinant does not change, since a determinant with two equal
columns vanishes \cite[Lemma~7]{Foata2007Han}, as in
Lemma~\ref{lem:Nprops}\,(b), and it is then expanded along its last
column \cite[Property~5]{Foata2007Han}, as in~\eqref{eq:Er}.  Here
the rows form an arbitrary set $L\cup u$, the last column $[X]$ is
arbitrary, and the expansion~\eqref{eq:El} along the first column is
needed as well.

\begin{proposition}\label{prop:Npath}
For every $L\subseteq J'$ and every letter~$\beta$, with $u$ fixed as
above,
\begin{equation}\label{eq:Npath}
N_L(\beta)\;\equiv\;\sum_{(v_1,\dots,v_s)}(-1)^{s}\;
\bw{v_s}{\beta}\bw{v_{s-1}}{v_s}\cdots\bw{v_1}{v_2}\bw{u}{v_1}\;
\det\bigl(\B_{L\setminus\{v_1,\dots,v_s\}}\bigr),
\end{equation}
where $(v_1,\dots,v_s)$, $s\ge0$, runs over the sequences of distinct
letters of~$L$ (for $s=0$ the word is $\bw{u}{\beta}$).
\end{proposition}

\begin{proof}
We argue by induction on~$|L|$, for all letters~$\beta$ at once.  For
$L=\varnothing$ we have $N_\varnothing(\beta)=\bw{u}{\beta}$, which is
the right-hand side of~\eqref{eq:Npath}, reduced to its term $s=0$.
Let $L\neq\varnothing$, and assume~\eqref{eq:Npath} for every proper
subset of~$L$ and every letter.  By~\eqref{eq:El} with $X=\{\beta\}$,
\[
N_L(\beta)\equiv\bw{u}{\beta}\det(\B_L)
-\sum_{v\in L}\bw{v}{\beta}N_{L\setminus v}(v).
\]
The first term $\bw{u}{\beta}\det(\B_L)$ is the term $s=0$ in the
right-hand side of~\eqref{eq:Npath}.  For $v\in L$, the induction
hypothesis for the set $L\setminus v$ and the letter~$v$ in place
of~$\beta$ reads
\[
N_{L\setminus v}(v)\;\equiv\;\sum_{(v_1,\dots,v_{s-1})}(-1)^{s-1}\,
\bw{v_{s-1}}{v}\cdots\bw{v_1}{v_2}\bw{u}{v_1}\,
\det\bigl(\B_{L\setminus\{v_1,\dots,v_{s-1},v\}}\bigr),
\]
where $(v_1,\dots,v_{s-1})$ runs over the sequences of distinct
letters of $L\setminus v$ (for $s=1$ the word is $\bw{u}{v}$).
Multiplying on the left by $-\bw{v}{\beta}$ and putting $v_s=v$, we
obtain exactly the terms of the right-hand side of~\eqref{eq:Npath}
with $s\ge1$ and last letter $v_s=v$.  Summing over $v\in L$ gives all
the terms with $s\ge1$.
\end{proof}

The next proposition is the heart of the proof.  Recall that
$\Delta_J$ is the forest sum defined in~\eqref{eq:Delta}.  The row
sum~$S_u$, which stands on the right in~$\Delta_J$, is absorbed into
the determinants as an additional column.

\begin{proposition}[absorption]\label{prop:absorb}
With $u$ fixed as above, we have
\[
\Delta_J\;\equiv\;\sum_{\beta\in\setA\setminus J}\ \sum_{L\subseteq J'}
(-1)^{|L|}\,N_L(\beta)\,S_{J'\setminus L}\,.
\]
\end{proposition}

\begin{proof}
Split the sum \eqref{eq:Delta} according to whether $u\in L$.  As the
row sums commute, $S_{J\setminus L}\equiv S_uS_{J'\setminus L}$ for
$L\subseteq J'$, and Lemma~\ref{lem:Nprops}\,(c) gives
\begin{equation}\label{eq:absplit}
\Delta_J\;\equiv\;\sum_{L\subseteq J'}(-1)^{|L|}
\Bigl(\det(\B_L)\,S_u-N_L(u)\Bigr)S_{J'\setminus L}\,.
\end{equation}
Fix $L\subseteq J'$ and compute $N_L(\setA)$ in two ways.
By~\eqref{eq:Er},
$N_L(\setA)\equiv\det(\B_L)S_u-\sum_{x\in L}N_{L\setminus x}(x)S_x$.
By Lemma~\ref{lem:Nprops}\,(a),(b), and since
$\setA=\{u\}\sqcup(\setA\setminus J)\sqcup(J'\setminus L)\sqcup L$,
$N_L(\setA)\equiv N_L(u)+\sum_{\beta\notin J}N_L(\beta)
+\sum_{x\in J'\setminus L}N_L(x)$.  Therefore
\[
\det(\B_L)\,S_u-N_L(u)\;\equiv\;
\sum_{\beta\in\setA\setminus J}N_L(\beta)
\;+\sum_{x\in J'\setminus L}N_L(x)
\;+\sum_{x\in L}N_{L\setminus x}(x)\,S_x.
\]
Insert the above relation into~\eqref{eq:absplit}.  The last two sums
contribute
\[
\sum_{L\subseteq J'}\ \sum_{x\in J'\setminus L}(-1)^{|L|}
N_L(x)\,S_{J'\setminus L}
\;+\;
\sum_{L\subseteq J'}\ \sum_{x\in L}(-1)^{|L|}
N_{L\setminus x}(x)\,S_x\,S_{J'\setminus L}\,.
\]
Writing $L=L_0\cup x$ in the second double sum and using
$S_xS_{J'\setminus L}\equiv S_{J'\setminus L_0}$, its term is
congruent to the opposite of the term $(L_0,x)$ of the first double
sum.  The two double sums cancel.
\end{proof}

\begin{proof}[Proof of Theorem~\textup{\ref{thm:pathintro}}]
Since the letter $u\in J$ fixed at the beginning of this section is
arbitrary, it suffices to prove the theorem for this~$u$.  By
Proposition~\ref{prop:absorb},
\[
\Delta_J\;\equiv\;\sum_{\beta\in\setA\setminus J}\ \sum_{L\subseteq J'}
(-1)^{|L|}\,N_L(\beta)\,S_{J'\setminus L}.
\]
Replacing each $N_L(\beta)$ by its expansion~\eqref{eq:Npath} from
Proposition~\ref{prop:Npath}, we get
\[
\Delta_J\;\equiv\;\sum_{\beta\in\setA\setminus J}\Xi(\beta),
\]
with
\[
\Xi(\beta)\;=\;\sum_{L\subseteq J'}\ \sum_{(v_1,\dots,v_s)}
(-1)^{|L|+s}\,w(P,\beta)\,\det(\B_{L\setminus V})\,S_{J'\setminus L},
\]
where $(v_1,\dots,v_s)$ runs over the sequences of distinct letters
of~$L$, $P=(u,v_1,\dots,v_s)$ and $V=\{u,v_1,\dots,v_s\}$; indeed
$L\setminus\{v_1,\dots,v_s\}=L\setminus V$ since $u\notin L$.  Exchange
the two sums in $\Xi(\beta)$: first choose the sequence
$(v_1,\dots,v_s)$ of distinct letters of~$J'$, that is, the escaping
path $(P,\beta)$ from~$u$ in~$J$, and then $L=L_0\cup\{v_1,\dots,v_s\}$
with $L_0\subseteq J\setminus V$.  Then
$J'\setminus L=(J\setminus V)\setminus L_0$ and
$(-1)^{|L|+s}=(-1)^{|L_0|}$, so that, by~\eqref{eq:Delta},
\begin{align*}
\Xi(\beta)&\;=\;\sum_{P}w(P,\beta)\sum_{L_0\subseteq J\setminus V}
(-1)^{|L_0|}\det(\B_{L_0})\,S_{(J\setminus V)\setminus L_0}\\
&\;=\;\sum_{P}w(P,\beta)\,\Delta_{J\setminus V}.
\end{align*}
Summing over $\beta\in\setA\setminus J$ gives the theorem.
\end{proof}

\section{Proof of the main theorem}\label{sec:proof}

\begin{proposition}[M\"obius inversion]\label{prop:mobius}
Recall that $\Delta_J$ is the forest sum defined in~\eqref{eq:Delta}.
For every finite $K\subseteq\setA$ we have
\begin{equation}\label{eq:rowsum}
\det(\B_K)\;\equiv\;\sum_{J\subseteq K}(-1)^{|J|}\,\Delta_J\,
S_{K\setminus J}.
\end{equation}
Conversely, if elements $\widetilde\Delta_J$ satisfy~\eqref{eq:rowsum}
with $\widetilde\Delta_J$ in place of~$\Delta_J$, for every finite~$K$, then $\widetilde\Delta_J\equiv\Delta_J$
for every~$J$.
\end{proposition}

\begin{proof}
Inserting \eqref{eq:Delta} into the right-hand side
of~\eqref{eq:rowsum} and using
$S_{J\setminus L}S_{K\setminus J}\equiv S_{K\setminus L}$ gives
$\sum_{L\subseteq K}(-1)^{|L|}\det(\B_L)S_{K\setminus L}
\sum_{L\subseteq J\subseteq K}(-1)^{|J|}$, and the inner sum is
$(-1)^{|K|}$ if $L=K$ and $0$ otherwise.  Only the term $L=K$
remains, and it is $(-1)^{|K|}\det(\B_K)(-1)^{|K|}=\det(\B_K)$.  This
proves~\eqref{eq:rowsum}.

For the converse, assume that
\[
\det(\B_L)\;\equiv\;\sum_{J\subseteq L}(-1)^{|J|}\,\widetilde\Delta_J\,
S_{L\setminus J}
\]
for every finite~$L$.  Insert this into~\eqref{eq:Delta}, written
for~$K$.  Since $S_{L\setminus J}S_{K\setminus L}\equiv
S_{K\setminus J}$, we get
\[
\Delta_K\;\equiv\;\sum_{J\subseteq K}(-1)^{|J|}\,\widetilde\Delta_J\,
S_{K\setminus J}\sum_{J\subseteq L\subseteq K}(-1)^{|L|}.
\]
The inner sum is $(-1)^{|K|}$ if $J=K$ and $0$ otherwise, so only the
term $J=K$ remains, and $\Delta_K\equiv\widetilde\Delta_K$.
\end{proof}

\begin{proof}[Proof of Theorem~\textup{\ref{thm:main}}]
We prove $\Delta_K\equiv\Gamma_K$ by induction on~$|K|$, the case
$K=\varnothing$ being trivial.  Let $K\neq\varnothing$, and assume
$\Delta_G\equiv\Gamma_G$ for every~$G$ with $|G|<|K|$.  For $v\in K$
put
\[
R_v\;=\;\sum_{(P,\beta)}w(P,\beta)\,\Gamma_{K\setminus V},
\]
where $(P,\beta)$ runs over the escaping paths from~$v$ in~$K$, of
support~$V$; all the biwords of~$R_v$ have top content~$K$.  By the
induction hypothesis and Theorem~\ref{thm:pathintro} for $J=K$ and
the letter~$v$,
\begin{equation}\label{eq:Rv}
R_v\;\equiv\;\sum_{(P,\beta)}w(P,\beta)\,\Delta_{K\setminus V}
\;\equiv\;\Delta_K\qquad(v\in K).
\end{equation}
By Lemma~\ref{lem:slices}\,(b), $(R_v)^{\langle v\rangle}
=\sum_{(P,\beta)}w(P,\beta)\,(\Gamma_{K\setminus V})^{\langle v\rangle}$,
so that Proposition~\ref{prop:firstchain} reads
$\Gamma_K=\sum_{v\in K}(R_v)^{\langle v\rangle}$.

Fix $u\in K$.  By~\eqref{eq:Rv}, $R_v\equiv R_u$ for every $v\in K$,
and slices preserve congruences, so that
\[
\Gamma_K\;\equiv\;\sum_{v\in K}(R_u)^{\langle v\rangle}
\;=\;R_u-\pi_*(R_u)
\]
by Lemma~\ref{lem:slices}\,(a).  Every biword of~$R_u$ has $|K|$
bottom letters, one of which is $\beta\notin K$, so it does not
cover~$K$ and $\pi_*(R_u)=0$.  Hence $\Gamma_K\equiv R_u\equiv\Delta_K$
by~\eqref{eq:Rv}.  The second statement of Theorem~\ref{thm:main}
follows by Proposition~\ref{prop:mobius}.
\end{proof}

\section{The quantum MacMahon Master Theorem}\label{sec:MMT}

Let $\widehat\calA$ be the algebra of formal series in the biletters.
Since $\calI$ is generated by homogeneous elements, congruences extend
to~$\widehat\calA$, a series being $\equiv0$ when all its homogeneous
components lie in~$\calI$.  Let $U_i=(1-S_i)^{-1}=\sum_{n\ge0}S_i^n$,
and $U_X$ the product of the $U_i$, $i\in X$.  Put
\[
\Ferm=\det(\I-\B)=\sum_{\sigma}(-1)^{\inv\sigma}
\prod_{j=1}^{r}\Bigl(\delta_{\sigma_j,j}-\bw{\sigma_j}{j}\Bigr),
\qquad
\Bos=\sum_{w}\bw{\overline w}{w},
\]
where $\sigma$ runs over the permutations of~$\setA$, the product is
taken in increasing order of~$j$, $w$ runs over all words on~$\setA$,
and $\overline w$ is the nondecreasing rearrangement of~$w$.

\begin{corollary}[quantum MacMahon Master Theorem \cite{Garoufalidis2006Le}]
\label{cor:MMT}
$\Ferm\cdot\Bos\equiv1$.
\end{corollary}

The proof goes through a formula for the product of $\Ferm$ with the
sum $U_\setA$ of all biwords with nondecreasing top word.

\begin{lemma}\label{lem:phi}
Let $\varphi$ be the algebra endomorphism of~$\calA$ with
$\varphi\bw{x}{a}=\delta_{xa}-\bw{x}{a}$.  Then
$\varphi(\calI)\subseteq\calI$.
\end{lemma}

\begin{proof}
It suffices to show that $\varphi$ fixes each generator
\eqref{R1}--\eqref{R2}.  Write $e_{xa}=\delta_{xa}$.  The image of the
generator \eqref{R2} is the generator itself (the quadratic part),
plus the linear part
\[
-\Bigl(e_{xa}\bw{y}{b}+e_{yb}\bw{x}{a}+e_{xb}\bw{y}{a}+e_{ya}\bw{x}{b}
-e_{ya}\bw{x}{b}-e_{xb}\bw{y}{a}-e_{yb}\bw{x}{a}-e_{xa}\bw{y}{b}\Bigr)
=0,
\]
plus the constant $e_{xa}e_{yb}+e_{xb}e_{ya}-e_{ya}e_{xb}-e_{yb}e_{xa}=0$.
The generator \eqref{R1} is treated in the same way.
\end{proof}

\begin{theorem}\label{thm:FU}
We have
\[
\Ferm\cdot U_\setA\;\equiv\;\sum_{J\subseteq\setA}\Gamma_J\,U_J.
\]
\end{theorem}

\begin{proof}
Apply~$\varphi$ to Theorem~\ref{thm:main} with $K=\setA$; this is
legitimate by Lemma~\ref{lem:phi}.  By definition
$\varphi(\det\B_\setA)=\Ferm$, and $\varphi(S_i)=1-S_i$.  An acyclic
map has no fixed point, so $\varphi(\Reg(b))=(-1)^{|J|}\Reg(b)$ for
every acyclic map~$b$ on~$J$, and $\varphi(\Gamma_J)=(-1)^{|J|}\Gamma_J$.  Hence
\[
\Ferm\;\equiv\;\sum_{J\subseteq\setA}\Gamma_J\,
\prod_{i\in\setA\setminus J}(1-S_i).
\]
Multiply on the right by $U_\setA$.  The row sums commute
modulo~$\calI$, hence so do all series in them, and $(1-S_i)U_i=1$;
therefore $\prod_{i\in\setA\setminus J}(1-S_i)\cdot U_\setA\equiv U_J$.
Substituting this in the previous congruence gives
$\Ferm\cdot U_\setA\equiv\sum_{J}\Gamma_J\,U_J$.
\end{proof}

\begin{proof}[Proof of Corollary~\textup{\ref{cor:MMT}}]
Call a biword \emph{real} if its top and bottom contents are equal,
and let $\Re$ be the projection on the real biwords; it preserves
congruences by~\eqref{eq:piI}.  Every biword of~$\Ferm$ is real, so
$\Re(\Ferm\cdot X)=\Ferm\cdot\Re(X)$ for every series~$X$, and
$\Re(U_\setA)=\Bos$, since $U_\setA$ is the sum of all biwords with
nondecreasing top word.  For $J\neq\varnothing$, every biword of
$\Gamma_JU_J$ has all its top letters in~$J$, while its $\Gamma_J$-part, an
acyclic map on~$J$, has a value outside~$J$; hence $\Re(\Gamma_JU_J)=0$.
Since $\Gamma_\varnothing=1$, taking real parts in
Theorem~\ref{thm:FU} gives $\Ferm\cdot\Bos\equiv1$.
\end{proof}

This is the right-quantum Master Theorem of~\cite{Garoufalidis2006Le} at $q=1$; the
general case follows by the ``$1=q$'' principle~\cite{Foata2008HanBasis}.
A bijective proof is given in \cite{Konvalinka2007Pak}; see also \cite{Konvalinka2008}.
In the commutative
case the same argument derives MacMahon's Master Theorem
\cite{MacMahon1960} from Theorem~\ref{thm:forest}, the forest version
of the matrix-tree theorem; we have not found this derivation in the
literature.

Theorem~\ref{thm:FU} describes the whole product
$\Ferm\cdot U_\setA$, and not only its real part: it expresses the
product of $\Ferm$ with all of $U_\setA$ as a sum with coefficients~$+1$, indexed by the
pairs formed by an acyclic map on a set~$J$ and a biword with
nondecreasing top word in the letters of~$J$.

\section{Cycle expansions}\label{sec:cycles}

The main theorem is now proved.  In this section we show that its
ideas extend from forests to maps with cycles.  When the escaping path
from a letter returns to that letter, it closes into a cycle; this
gives an expansion of the determinant at an arbitrary letter
(Section~\ref{ssec:atletter}), a recursion for the forest sums
(Section~\ref{ssec:recursion}), a cycle notation formula for
$\det(\B_M)$ (Section~\ref{ssec:cyclenot}), and an identity that
interpolates between Theorem~\ref{thm:main} and this formula
(Section~\ref{ssec:interp}).

\subsection{Expansion at a letter}\label{ssec:atletter}
For $\beta=u$, the paths of Proposition~\ref{prop:Npath} close into
cycles.  For a cycle $Z:\ u\to v_1\to\dots\to v_s\to u$ of distinct
letters put
\begin{equation}\label{eq:rho}
\rho_u(Z)=\bw{v_s}{u}\bw{v_{s-1}}{v_s}\cdots\bw{v_1}{v_2}\bw{u}{v_1}
\qquad(\rho_u(Z)=\bw{u}{u}\ \text{if}\ s=0).
\end{equation}
Applying Proposition~\ref{prop:Npath} to $J=K\cup u$, $L=K$ and
$\beta=u$, where $K\subseteq\setA$ is finite and $u\notin K$, and
using $N_K(u)\equiv\det(\B_{K\cup u})$ (Lemma~\ref{lem:Nprops}\,(c))
and $(-1)^s=(-1)^{|Z|-1}$, we obtain:

\begin{corollary}[cycle expansion]\label{cor:cyclepeel}
Let $K\subseteq\setA$ be finite and $u\in\setA\setminus K$.  Then
\[
\det\bigl(\B_{K\cup u}\bigr)
\;\equiv\;
\sum_{Z}\,(-1)^{|Z|-1}\,\rho_u(Z)\,\det\bigl(\B_{K\setminus Z}\bigr),
\]
where $Z$ runs over the cycles through~$u$ with support in $K\cup u$.
\end{corollary}

The resummation of the proof of Theorem~\ref{thm:pathintro}, applied
to~\eqref{eq:absplit} and Corollary~\ref{cor:cyclepeel}, gives the
following recursion; here $J\subseteq\setA$ is finite and, for
$u\in J$, $J'=J\setminus u$, as in Section~\ref{sec:path}.

\begin{proposition}[cycle recursion]\label{prop:cyclerecu}
For every $u\in J$,
\[
\Delta_J\;\equiv\;\Delta_{J\setminus u}\,S_u
\;-\;\sum_{Z\ni u}\rho_u(Z)\,\Delta_{J\setminus Z}\,,
\]
where $Z$ runs over the cycles through~$u$ with support in~$J$.
\end{proposition}

\begin{proof}
By \eqref{eq:absplit} and Lemma~\ref{lem:Nprops}\,(c),
\[
\Delta_J\equiv\Delta_{J'}S_u
-\sum_{L\subseteq J'}(-1)^{|L|}\det(\B_{L\cup u})S_{J'\setminus L}.
\]
Expand $\det(\B_{L\cup u})$ by Corollary~\ref{cor:cyclepeel} and put
$L_0=L\setminus Z$; the signs combine to $(-1)^{|L_0|}$, and the sum
over $L_0\subseteq J\setminus Z$ is $\rho_u(Z)\Delta_{J\setminus Z}$.
\end{proof}

\subsection{A recursion for the forest sums}\label{ssec:recursion}
When $u=\max Z$, we write $\rho(Z)=\rho_u(Z)$: the cycle written
backwards from its maximum.  For an acyclic map, a new letter~$x$
larger than all others can be attached by multiplying by~$S_x$ on the
right; the result is a sum of acyclic maps and of maps with a single
cycle, passing through~$x$, each written as its normal word \eqref{eq:W}.

\begin{corollary}\label{cor:RL}
Let $J\subseteq\setA$ be finite and $x\in\setA$, $x>\max J$.  Then
\[
\Gamma_J\,S_x\;\equiv\;\Gamma_{J\cup x}\;+\;\sum_{Z\ni x}\rho(Z)\,
\Gamma_{(J\cup x)\setminus Z}\,,
\]
the sum being over the cycles through~$x$ with support in $J\cup x$.
\end{corollary}

\begin{proof}
Apply Proposition~\ref{prop:cyclerecu} to $J\cup x$ and $u=x$.  Since
$x$ is the largest letter of every cycle~$Z$ through~$x$, we have
$\rho_x(Z)=\rho(Z)$, and the proposition reads
$\Delta_J S_x\equiv\Delta_{J\cup x}
+\sum_{Z\ni x}\rho(Z)\,\Delta_{(J\cup x)\setminus Z}$.  By
Theorem~\ref{thm:main}, $\Delta_G\equiv \Gamma_G$ for every~$G$, and
replacing each~$\Delta$ by the corresponding~$\Gamma$ gives the corollary.
\end{proof}

Since a map on $J\cup x$ with a unique cycle $Z\ni x$ is the same as
the cycle together with an acyclic map on $(J\cup x)\setminus Z$, the
right-hand side is the sum, over the maps on $J\cup x$ all of whose
cycles pass through~$x$, of the cycle word followed by the regular
form of the rest.

\subsection{The cycle notation expansion}\label{ssec:cyclenot}
Iterating Corollary~\ref{cor:cyclepeel} at the largest letter
expands the determinant over cycle decompositions.  For a permutation
$\sigma$ of a finite set~$M$ list its cycles $Z_1,\dots,Z_c$ in
decreasing order of their maxima.  Since
$\sum_i(|Z_i|-1)=|M|-c$:

\begin{corollary}\label{cor:cycleform}
For every finite $M\subseteq\setA$,
\[
\det(\B_M)\;\equiv\;(-1)^{|M|}
\sum_{\sigma}\,(-1)^{\cyc\sigma}\,
\rho(Z_1)\,\rho(Z_2)\cdots\rho(Z_c),
\]
where $\sigma$ runs over the permutations of~$M$, with cycles
$Z_1,\dots,Z_c$ listed in decreasing order of their maxima.
\end{corollary}

For instance, for $M=\{1,2\}$ the right-hand side is
$\bw{21}{21}-\bw{12}{21}$, and \eqref{R2} turns it into
$\det(\B_M)=\bw{12}{12}-\bw{21}{12}$.

\begin{remark}\label{rem:Konv}
By Lemma~\ref{lem:rev}, Corollary~\ref{cor:cycleform} is equivalent to
\[
\det(\B_M)\;\equiv\;(-1)^{|M|}
\sum_{\sigma}\,(-1)^{\cyc\sigma}\,
\omega(\rho(Z_c))\cdots\omega(\rho(Z_2))\,\omega(\rho(Z_1)),
\]
where $\omega(\rho(Z))=\bw{u}{v_1}\bw{v_1}{v_2}\cdots\bw{v_s}{u}$ for
the cycle $Z:\ u\to v_1\to\dots\to v_s\to u$ with $u=\max Z$.  In the
language of lattice paths of \cite{Konvalinka2007Pak,Konvalinka2008}, $\omega(\rho(Z))$ is
the closed lattice path that starts at the maximum of~$Z$ and goes
once around~$Z$, and the cycles are listed in increasing order of
their maxima.  In this form the identity is essentially Konvalinka's
expansion of the determinant of a right-quantum matrix into
(back-)path sequences \cite[Section~4]{Konvalinka2008}; here it comes out of
the path expansion.  For $M=\{1,2\}$ the identity reads
$\det(\B_M)\equiv\bw{12}{12}-\bw{21}{12}$, which is the definition of
$\det(\B_M)$.
\end{remark}

\subsection{An interpolating identity}\label{ssec:interp}
Theorem~\ref{thm:main} and Corollary~\ref{cor:cycleform} are the two
extreme cases of one identity.  For a map~$b$ on a finite set, with
cycles $Z_1,\dots,Z_c$ listed in decreasing order of their maxima and
restriction $\hat b$ to the letters on no cycle, define the \emph{normal
word} of~$b$ by
\begin{equation}\label{eq:W}
W(b)\;=\;\rho(Z_1)\cdots\rho(Z_c)\;\Reg(\hat b).
\end{equation}

\begin{example}\label{ex:treecycles}
Extend the map of Example~\ref{ex:tree} to the set
$\{1,\dots,14\}\setminus\{11\}$ by $b(2)=2$, $b(4)=8$ and $b(8)=4$.
Its functional graph is drawn below: the forest of
Example~\ref{ex:tree}, whose root~$8$ now lies on the cycle
$4\rightleftarrows8$, together with the fixed point~$2$; only the
root~$11$ lies outside the domain.
\begin{center}
\begin{tikzpicture}[>=Stealth,scale=0.95,every node/.style={font=\small}]
\node (n11) at (0,0) [draw,rectangle] {$11$};
\node (n3) at (-0.9,1.1) {$3$};
\node (n9) at (0.9,1.1) {$9$};
\draw[->] (n3) -- (n11);
\draw[->] (n9) -- (n11);
\node (n8) at (4,0) {$8$};
\node (n4) at (4,-1.1) {$4$};
\draw[->] (n8) to[bend left=35] (n4);
\draw[->] (n4) to[bend left=35] (n8);
\node (n12) at (4,1.1) {$12$};
\draw[->] (n12) -- (n8);
\node (n14) at (3,2.2) {$14$};
\node (n6) at (5.3,2.2) {$6$};
\draw[->] (n14) -- (n12);
\draw[->] (n6) -- (n12);
\node (n1) at (2.2,3.3) {$1$};
\node (n5) at (3.8,3.3) {$5$};
\draw[->] (n1) -- (n14);
\draw[->] (n5) -- (n14);
\node (n13) at (5.3,3.3) {$13$};
\draw[->] (n13) -- (n6);
\node (n7) at (1.5,4.4) {$7$};
\node (n10) at (2.9,4.4) {$10$};
\draw[->] (n7) -- (n1);
\draw[->] (n10) -- (n1);
\node (n2) at (7,0) {$2$};
\draw[->] (n2) to[out=60,in=120,looseness=8] (n2);
\end{tikzpicture}
\end{center}
The cycles, in decreasing order of their maxima, are
$Z_1=\{4,8\}$ and $Z_2=\{2\}$, with $\rho(Z_1)=\bw{4\;8}{8\;4}$ and
$\rho(Z_2)=\bw{2}{2}$.  The letters on no cycle form the set~$J$ of
Example~\ref{ex:tree}, and $\hat b$ is the map of that example.
Hence
\[
W(b)\;=\;\underbrace{\bw{4\;8}{8\;4}\,\bw{2}{2}}_{\text{cycles}}\;
\underbrace{\bw{3}{11}\bwc{12/8,14/12,5/14}
\bwc{1/14,7/1}\bw{9}{11}\bw{10}{1}\bwc{6/12,13/6}}
_{\Reg(\hat b)}\,.
\]
This example only illustrates the normal word~$W(b)$: it is not one
of the maps of Theorem~\ref{thm:master} below, in which every cycle
must meet a set~$M$ of letters larger than all the letters
outside~$M$.
\end{example}

For disjoint finite sets $J,M\subseteq\setA$ with $\max J<\min M$
(when both are nonempty) put
\[
V_J(M)\;=\;\sum_{L\subseteq J}(-1)^{|L|}\,
\det\bigl(\B_{L\cup M}\bigr)\;S_{J\setminus L}\,.
\]

\begin{theorem}\label{thm:master}
With $J,M$ as above,
\[
V_J(M)\;\equiv\;(-1)^{|M|}\sum_{b}\;(-1)^{\cyc(b)}\;W(b),
\]
where $b$ runs over the maps on $J\cup M$ whose cycles all meet~$M$
and cover~$M$, and $\cyc(b)$ is the number of cycles of~$b$.
\end{theorem}

For $M=\varnothing$ this is Theorem~\ref{thm:main}, and for
$J=\varnothing$ it is Corollary~\ref{cor:cycleform}.  The identity
thus interpolates between the forest expansion, in which no order of
the vertices fixed in advance is admissible, and the cycle expansion
of Remark~\ref{rem:Konv}, in which such an order is.

\begin{proof}
Induction on $|J|$, the case $J=\varnothing$ being
Corollary~\ref{cor:cycleform}.  Let $y=\max J$ and $J'=J\setminus y$.
Splitting according to whether $y\in L$, and since $y$ is the largest
letter of $J\setminus L$ when $y\notin L$,
\begin{equation}\label{eq:Vrec}
V_J(M)\;=\;V_{J'}(M)\,S_y\;-\;V_{J'}(M\cup y),
\end{equation}
and both terms are covered by the induction hypothesis.  After
cancelling $(-1)^{|M|}$, the claim becomes
\begin{equation}\label{eq:step}
\sum_{b'}(-1)^{\cyc(b')}\,W(b')\,S_y
\;+\;\sum_{b\in\mathcal{M}_{M\cup y}}(-1)^{\cyc(b)}\,W(b)
\;\equiv\;
\sum_{b\in\mathcal{M}_M}(-1)^{\cyc(b)}\,W(b),
\end{equation}
where $b'$ runs over the maps on $J'\cup M$ whose cycles meet and
cover~$M$, and $\mathcal{M}_N$ denotes the set of maps on $J\cup M$
whose cycles meet and cover~$N$.  The pairs $(b',a)$, $a\in\setA$,
correspond to the maps $b=b'\cup(y\mapsto a)$ on $J\cup M$ whose
cycles avoiding~$y$ meet and cover~$M$.  Three cases occur.

(i) $y$ lies on no cycle of~$b$: then $b$ comes from a pair and lies
in $\mathcal{M}_M\setminus\mathcal{M}_{M\cup y}$, with
$\cyc(b)=\cyc(b')$.

(ii) $y$ lies on a cycle of~$b$ meeting~$M$: then $b$ does not come
from a pair, and $b\in\mathcal{M}_M$ if and only if
$b\in\mathcal{M}_{M\cup y}$; these terms are the same on both sides of
\eqref{eq:step}.

(iii) $y$ lies on a cycle $Z_y$ of~$b$ avoiding~$M$: then
$Z_y\subseteq J$ and $\max Z_y=y$; such a~$b$ lies in
$\mathcal{M}_{M\cup y}\setminus\mathcal{M}_M$ if and only if it comes
from a pair, and $\cyc(b)=\cyc(b')+1$.  As the other cycles of~$b$
have their maxima in~$M$, the factor $\rho(Z_y)$ stands last among the
cycle factors of~$W(b)$.

In cases (i) and (iii) all the words begin with the factor
$\rho(Z_1)\cdots\rho(Z_c)$ formed by the cycles of~$b'$, and the
remaining letters form a set $Q\subseteq J$ with $y=\max Q$.  Removing
this common left factor, \eqref{eq:step} becomes, for each choice of
$Z_1,\dots,Z_c$, the congruence of Corollary~\ref{cor:RL} for
$(Q\setminus y,\,y)$, which holds.  Together with case~(ii), this
proves \eqref{eq:step}, and hence the theorem.
\end{proof}

\section{Further consequences}\label{sec:cons}

\subsection{The bottom row and Cayley's formula}\label{ssec:cayley}
Let $\mathbb{Z}\langle y_1,\dots,y_r\rangle$ be the free associative
algebra on noncommuting variables $y_a$, $a\in\setA$, and let
$\theta:\calA\to\mathbb{Z}\langle y_1,\dots,y_r\rangle$ be the algebra
morphism with $\theta\bw{x}{a}=y_a$, which keeps only the bottom row
of a biword.  It maps \eqref{R1} to $y_ay_a-y_ay_a=0$ and \eqref{R2}
to $y_ay_b+y_by_a-y_ay_b-y_by_a=0$, so $\theta(\calI)=0$; in other
words, the matrix whose rows are all equal to $(y_1,\dots,y_r)$ is a
Manin matrix.  For $X\subseteq\setA$ put $Y_X=\sum_{a\in X}y_a$ and
$Y=Y_\setA$.

\begin{corollary}\label{cor:cayley}
Let $J\subseteq\setA$ be nonempty, $R=\setA\setminus J$ and $m=|J|$.
Then
\[
\sum_{b:J\to\setA\ \mathrm{acyclic}}\theta\bigl(\Reg(b)\bigr)
\;=\;Y_R\,Y^{m-1}
\qquad\text{in }\mathbb{Z}\langle y_1,\dots,y_r\rangle.
\]
Equivalently, the bottom row of $\Reg(b)$ is a bijection from the
acyclic maps on~$J$ onto the words $c_1c_2\cdots c_m$ on~$\setA$ with
$c_1\in R$.
\end{corollary}

\begin{proof}
We have $\theta(S_i)=Y$ for every~$i$.  For $|L|\ge2$, all the terms
of $\det(\B_L)$ have the same bottom word $c_1\cdots c_m$, and their
signs add up to zero, so $\theta(\det\B_L)=0$; moreover
$\theta(\det\B_{\{j\}})=y_j$.  Hence
$\theta(\Delta_J)=Y^m-Y_JY^{m-1}=Y_RY^{m-1}$, and the identity is the
image of Theorem~\ref{thm:main} under~$\theta$.  Its right-hand side
is the sum of the words $c_1\cdots c_m$ with $c_1\in R$, each with
coefficient one, and its left-hand side is a sum of words with
nonnegative coefficients, one for each acyclic map; this is the
second statement.
\end{proof}

For an acyclic map~$b$ on~$J$, let $f$ be the map of~$\setA$ that
extends~$b$ by fixing every letter of~$R$.  By Remark~\ref{rem:FF},
the Foata--Fuchs word $\phi(f)$ is the increasing word of~$R$ followed
by the bottom row of $\Reg(b)$.  Thus Corollary~\ref{cor:cayley} says
that $\phi$ maps these maps~$f$ bijectively onto the words formed by
the increasing word of~$R$ followed by a word $c_1\cdots c_m$ with
$c_1\in R$: this is the case of forests in
\cite[Theorem~5.9 and (6.1)]{Foata1970Fuchs}, obtained here as the
bottom-row image of Theorem~\ref{thm:main}.  This bijection for
forests with a given set of roots is also described in
\cite[Section~2]{AddarioBerry2023}, where the Foata--Fuchs bijection
is applied to random trees.  Letting the~$y_a$ commute gives the
multivariate Cayley formula
$\sum_b\prod_{v\in J}y_{b(v)}=Y_R\,Y^{m-1}$, and $y_a=1$ gives
$|R|\,r^{m-1}$ forests, hence $r^{r-2}$ trees on $r$ vertices
converging to a given root.  Thus Theorem~\ref{thm:main} contains
both the matrix-tree theorem, in its commutative image, and the
Foata--Fuchs encoding of forests, in its bottom-row image.

The same argument applies to every admissible ordering: if an
arrangement $\mathrm{Arr}(b)$ of each acyclic map~$b$ on~$J$ is chosen
so that $\Delta_J\equiv\sum_b\mathrm{Arr}(b)$, then the bottom row of
$\mathrm{Arr}(b)$ is a bijection onto the words $c_1\cdots c_m$ with
$c_1\in R$.  The orderings $\Ord_\varsigma$ and $\Reg_\pi$ of
Corollaries~\ref{cor:selector} and~\ref{cor:anyorder} below thus give
bijective proofs of Cayley's formula.  Conversely, this necessary
condition explains Example~\ref{ex:intro}: in $J=\{1,2\}$, $r=3$,
writing the arcs in the order $1,2$ (resp.\ $2,1$) gives the bottom
words $23,31,33$ (resp.\ $32,13,33$), and $23$ (resp.\ $13$) does not
begin with the root~$3$.

\subsection{Other admissible orderings}
A \emph{selector} is a rule~$\varsigma$ choosing a letter
$\varsigma(G)\in G$ for every nonempty finite $G\subseteq\setA$.  For
an acyclic map~$b$ on~$G$ define $\Ord_\varsigma(b)$ recursively:
it is empty if $G=\varnothing$; otherwise, with $u=\varsigma(G)$ and
$(P,\beta)$ the forward orbit of~$u$, of support~$V$,
$\Ord_\varsigma(b)=w(P,\beta)\,\Ord_\varsigma(b|_{G\setminus V})$.

\begin{corollary}\label{cor:selector}
For every selector~$\varsigma$ and every finite $J\subseteq\setA$,
\[
\Delta_J\;\equiv\;\sum_b\Ord_\varsigma(b),
\]
the sum being over the acyclic maps on~$J$.
\end{corollary}

\begin{proof}
Induction on~$|J|$.  With $u=\varsigma(J)$, the acyclic maps on~$J$
correspond bijectively to the pairs formed by an escaping path
$(P,\beta)$ from~$u$ in~$J$, of support~$V$, and an arbitrary acyclic
map on $J\setminus V$.  Hence $\sum_b\Ord_\varsigma(b)
=\sum_{(P,\beta)}w(P,\beta)\sum_{b'}\Ord_\varsigma(b')
\equiv\sum_{(P,\beta)}w(P,\beta)\Delta_{J\setminus V}\equiv\Delta_J$
by induction and Theorem~\ref{thm:pathintro}.
\end{proof}

For $\varsigma(G)=\min G$ one writes the escaping path of the smallest
letter, then the path of the smallest letter not yet written up to the
first letter already written, and so on.

For $\varsigma(G)=\max G$, $\Ord_\varsigma$ is the Pr\"ufer code
\cite{Prufer1918} read backwards.  For an acyclic map~$b$ on~$G$,
remove repeatedly the smallest leaf, that is, the smallest letter of
the remaining domain that is not a value of the remaining map, and
list the removed biletters $\bw{v}{b(v)}$; read backwards, this list
is $\Ord_\varsigma(b)$.  Indeed, let $u=\max G$ and let
$P=(u,b(u),\dots,b^s(u))$, $\beta=b^{s+1}(u)$, be the forward orbit
of~$u$, of support~$V$.  Since $b$ maps $V$ into $V\cup\{\beta\}$, no
letter of $G\setminus V$ has a preimage in~$V$, so the leaves in
$G\setminus V$ are those of $b|_{G\setminus V}$, and one of them exists
as long as letters of $G\setminus V$ remain, since $b$ is acyclic.
While no letter of~$V$ is removed, $u$ is the only possible leaf
in~$V$, and it is larger than every other letter.  Hence the letters
of $G\setminus V$ are removed first, in the order of the Pr\"ufer code
of $b|_{G\setminus V}$, and then $u,b(u),\dots,b^s(u)$, which gives
the claim by induction on~$|G|$.  For this ordering, the bottom-row
bijection of Section~\ref{ssec:cayley} is Pr\"ufer's bijection, read
backwards.

The regular form is not a selector ordering, since its first chain starts at the smallest source, which
depends on the map.  We have not been able to find other natural
admissible orderings, for instance by sorting the letters by their
depth in the forest, or by placing greedily the smallest (or the
largest) letter whose image is already placed.

The natural order of the alphabet plays no special role in
Theorem~\ref{thm:main}.  For a total order~$\pi$ of~$\setA$, let
$\Reg_\pi(b)$ be defined as in Definition~\ref{def:reg}, with~$\pi$ in
place of the natural order: each letter is assigned to its
$\pi$-smallest source, and the chains are listed in $\pi$-increasing
order of their sources.

\begin{corollary}\label{cor:anyorder}
For every total order~$\pi$ of~$\setA$ and every finite
$J\subseteq\setA$,
\[
\Delta_J\;\equiv\;\sum_b\Reg_\pi(b),
\]
the sum being over the acyclic maps on~$J$.
\end{corollary}

\begin{proof}
Let $\phi$ be the permutation of~$\setA$ such that $x$ precedes~$y$
in~$\pi$ if and only if $\phi(x)<\phi(y)$, and extend it to the
automorphism $\bw{x}{a}\mapsto\bw{\phi(x)}{\phi(a)}$ of~$\calA$.  It
maps the generators \eqref{R1}--\eqref{R2} to generators, so
$\phi(\calI)=\calI$, and $\phi(S_i)=S_{\phi(i)}$.  Let
$L=\{c_1<\dots<c_m\}$, $\phi(L)=\{c'_1<\dots<c'_m\}$ and
$\phi(c_k)=c'_{\lambda(k)}$, where $\lambda$ is a permutation of
$\{1,\dots,m\}$.  If $\sigma_k=c_{\alpha(k)}$ and
$\phi(\sigma_k)=c'_{\tau(k)}$, then $\tau=\lambda\alpha$, so that
$(-1)^{\inv\sigma}=\mathrm{sgn}(\alpha)
=\mathrm{sgn}(\lambda)\,\mathrm{sgn}(\tau)$.  Hence
\[
\phi\bigl(\det(\B_L)\bigr)=\mathrm{sgn}(\lambda)\,
\det\bigl(\phi(L);[c'_{\lambda(1)}],\dots,[c'_{\lambda(m)}]\bigr)
\equiv\det\bigl(\B_{\phi(L)}\bigr)
\]
by Lemma~\ref{lem:gendet}\,(ii), and therefore
$\phi(\Delta_J)\equiv\Delta_{\phi(J)}$.  The map $b\mapsto\phi b\phi^{-1}$
is a bijection from the acyclic maps on~$J$ onto the acyclic maps
on~$\phi(J)$, and $\phi^{-1}\bigl(\Reg(\phi b\phi^{-1})\bigr)=\Reg_\pi(b)$,
since $\phi^{-1}$ carries the natural order to~$\pi$.  Applying
$\phi^{-1}$ to Theorem~\ref{thm:main} for $\phi(J)$ gives the result.
\end{proof}

Thus $\Reg_\pi(b)$ is the Foata--Fuchs encoding with respect to the
order~$\pi$.  For the reverse order, each letter is assigned to its
largest source and the chains are listed in decreasing order of their
sources.  For
every~$\pi$, the order of the arcs in $\Reg_\pi(b)$ still depends on
the shape of the forest.

\begin{example}\label{ex:orders}
Let $r=4$ and $J=\{1,2,3\}$, and compare $\Reg$ with $\Ord_\varsigma$
for $\varsigma(G)=\max G$, that is, with the Pr\"ufer code read
backwards.  Only two acyclic maps on~$J$ have bottom
content $\{2,4,4\}$: $(b(1),b(2),b(3))=(2,4,4)$, the forest
$1\to2\to4$, $3\to4$, and $(4,4,2)$, the forest $1\to4$,
$3\to2\to4$.  Their words are
\[
\Reg:\ \bw{213}{424},\ \bw{123}{442};\qquad
\Ord_\varsigma:\ \bw{321}{442},\ \bw{231}{424}.
\]
We check that
$\bw{213}{424}+\bw{123}{442}\equiv\bw{321}{442}+\bw{231}{424}$.
\end{example}

\section*{Acknowledgments}
The author thanks Fr\'ed\'eric Chapoton, Dominique Foata, Kathy Q. Ji and
Huan Xiong for helpful discussions, and especially Dominique Foata for pointing out the reference \cite{Foata1970Fuchs}.

\section*{Use of AI}
The main result of this paper, Theorem~\ref{thm:main}, with the
regular form, was observed and conjectured by the author before 2017.  The AI assistant Claude
(Anthropic) was used, under the author's direction, to explore proof
strategies and to edit the text.

\bibliographystyle{amsplain}
\providecommand{\bysame}{\leavevmode\hbox to3em{\hrulefill}\thinspace}
\providecommand{\MR}{\relax\ifhmode\unskip\space\fi MR }
\providecommand{\MRhref}[2]{\href{http://www.ams.org/mathscinet-getitem?mr=#1}{#2}}

\end{document}